\magnification=\magstep1
\baselineskip=16pt  
\parindent=12pt
 \parskip=2pt 
\def\sqr#1#2{{\vcenter{\vbox{\hrule height.#2pt
     \hbox{\vrule width.#2pt height#1pt \kern#1pt
   \vrule width.#2pt}\hrule height.#2pt}}}}
\def\qed{\mathchoice\sqr64\sqr64\sqr{2.1}3\sqr{1.5}3} 
\def\ZZ{Z\!\!\!Z\,}
\def\RR{R\!\!\!\!\!I\,\,}

\def\11{l\!\!\!1\,}

\font\tenib=cmmib10
\newfam\mitbfam
\textfont\mitbfam=\tenib
\scriptfont\mitbfam=\seveni
\scriptscriptfont\mitbfam=\fivei



\def\and{ \hbox{ and } }

\def\a{\alpha}

\def\[{\left [ \  }
\def\]{\ \right ]  }
\def\){\  \right ) }
\def\({ \left (  \  }

\def\to{\rightarrow}

\def\frac{\over}
\def\\{\cr}
\def\ref{}

\hbox{}
\overfullrule=0in

\def \a {\alpha}

\def\phi{\varphi}
%
\hbox{}
\centerline{\bf Comparison inequality and two block estimate}
\centerline{\bf for inhomogeneous Bernoulli measures}\smallskip
\smallskip
\centerline{Jeremy Quastel}
\centerline{University of Toronto}
\bigskip
{\bf Abstract.}  We consider
inhomogeneous Bernoulli measures of the form $\prod_{x\in\Lambda} p_x$
where $p_x$ are prescribed and uniformly bounded above and below away from $0$ and
$1$.  A comparison inequality is proved between the Kawasaki and Bernoulli-Laplace Dirichlet forms.
Together with a recent result of Caputo on the gap of the Bernoulli-Laplace model, this proves
a spectral gap of the correct order $L^{-2}$ on cubes of side length $L$ for the Kawasaki dynamics.
The two block estimate of hydrodynamic limits is also obtained.

\bigskip
\noindent{\bf  0.  Introduction}
\medskip
Recently there has been a lot of interest in the transport properties of particle systems in 
random media [AHL, BE, F, GP, KPW, K, MA, R, Se].  
A simple model for which the hydrodynamic scaling limit can be obtained
is the Kawasaki dynamics for random Bernoulli measures studied in [Q], [QY] and recently, [FM].  Such systems have been used to
model electron transport in doped crystals.  The hydrodynamic limit is a nonlinear diffusion
equation with a nontrivial density dependent diffusion coefficient given by a Green-Kubo formula.
In fact for such a system even the existence of a diffusive scaling limit is nontrivial.  The
key input is a so called moving particles lemma which leads to a diffusive spectral gap, as well as a two block estimate [GPV].  
The hydrodynamic limit was described in [Q] and an unpublished manuscript [QY] contains most of the details of the proof.  
Based on these [FM] recently give a complete proof.  However they
do not use the long jumps method sketched in [Q], [QY], but the more traditional non-gradient method
together with subtraction of an appropriate term to compensate for the inhomogeneity in the medium.
This restricts one to dimensions $d\ge 3$.  Another recent article [C] proves the gap of the 
Bernoulli-Laplace model.  Together with a comparison inequality from [QY] he concludes the
spectral gap of the Kawasaki dynamics. 
The purpose of this note is to complete the literature by providing the needed moving particles lemma from the
unpublished manuscript [QY].
This work as well as [Q] and [QY] on the
hydrodynamic limit, arose out of problems suggested by Herbert Spohn and from joint work with
H. T. Yau.  Their contribution is
gratefully acknowledged.
 \bigskip
\noindent{\bf 1. Inhomogeneous Bernoulli measures.}
\medskip
Let $\a_x\in [-K,K]$, $x\in \ZZ^d$ be given and let
$$
p_x= {e^{\a_x}\over 1+ e^{\a_x}}, \qquad x\in \ZZ^d.
$$
For any $\Lambda\subset\ZZ^d$ the inhomogeneous Bernoulli
measure
$
\mu_\Lambda(\eta)$ on $\{0,1\}^\Lambda$ is given by
$$
\mu_\Lambda(\eta) = \prod_{x\in\Lambda} p_x^{\eta_x} (1-p_x)^{1-\eta_x}
=Z_{\a,\Lambda}^{-1}\exp\{- H_{\alpha,\Lambda}(\eta)\}, \qquad
H_{\alpha,\Lambda}(\eta)=-\sum_{x\in\Lambda}\a_x \eta_x.
$$
$Z_{\a,\Lambda}=\prod_{x\in\Lambda} (1+e^{\alpha_x})$ is the normalization to make it a probability measure. 
We can also condition to have a fixed number $N$ of particles in $\Lambda$,
$$
\mu_{\a,\Lambda,N}(\eta) = \mu_{\a,\Lambda}(\eta ~|~ \sum_x \eta_x =N ).
$$
For each configuration $\eta$ with exactly $N$ particles, $
\mu_{\a,\Lambda,N}(\eta) =  Z_{\a,\Lambda,N}^{-1}\exp\{- H_{\alpha,\Lambda}
(\eta)\} 
$
with
$$
 Z_{\a,\Lambda,N} = \sum_{A\subset\Lambda, |A|=N }\exp
 \{\sum_{x\in A}\alpha_x\}= {|\Lambda|\choose N}\zeta_{\a,\Lambda,N}, \qquad
\zeta_{\a,\Lambda,N}= Av_{A\subset\Lambda, |A|=N }\exp
 \{\sum_{x\in A}\alpha_x\} 
.
$$
Note that  
$$
\mu_{\alpha, \Lambda} (\eta_x=1) = {e^{\alpha_x}\over 1+e^{\alpha_x}}= p_x
$$
while the corresponding quantity in the canonical ensemble is
$$
\mu_{\alpha, \Lambda, N} (\eta_x=1) =e^{\a_x} {Z_{\a,\Lambda\setminus \{x\}, N-1}\over Z_{\a,\Lambda,N} }
={N\over |\Lambda|}\left( e^{\a_x}{\zeta_{\a,\Lambda\setminus \{x\}, N-1}\over \zeta_{\a,\Lambda,N} }\right)
 = p_{\alpha, \Lambda, N, x}.
$$
We will use the notation
$E_{\mu}[ f; g]$ to denote the covariance $E_\mu [ (f-E_\mu[f])(g-E_\mu[g])]$
as well as ${\rm Var}_\mu(f)$ for the variance $E_{\mu}[ f;f]$.

\bigskip
\noindent{\bf 2. Dynamics and main results}
\medskip

We define dynamics through Dirichlet forms.  There are two basic
dynamics: Glauber and Kawasaki.

{\bf Glauber.}  Let $\Lambda\subset \ZZ^d$.  The Dirichlet form is given
by
$$
D_{G}(\Lambda; f)= E_{\mu_\Lambda}
[ \sum_{x\in\Lambda}  (f(\sigma_x \eta) - f(\eta))^2]
$$
where
$$
(\sigma_x\eta)_y = \cases{ \eta_y & if $ y\neq x$; \cr 1-\eta_y & if $x=y$.}
$$
The dynamics corresponding to this Dirichlet form is when each site changes
its value from $\eta_x$ to $1-\eta_x$ at rate $ 1+ \exp\{ \a_x (1
-2\eta_x)\}$.

{\bf Kawasaki.}  On a connected set $\Lambda\subset\ZZ^d$ (in the
sense of nearest neighbours, which we write as $x\sim y$),
and given a fixed number of particles $0\le N\le |\Lambda|$, the Dirichlet
form is
$$
D_{Kaw}(\Lambda, N; f)=E_{\mu_{\Lambda,N}}
[ \sum_{x\sim y\atop
x,y\in\Lambda} (f(T_{x,y} \eta) - f(\eta))^2]
$$
where
$$ 
(T_{x,y} \eta )_z = 
\cases { \eta_y, & if $z=x$ \cr 
                                 \eta_x, & if $z=y$ \cr
                                  \eta_z, & otherwise. }
$$
 The corresponding dynamics is a system of $N$ particles on 
$\Lambda=\Lambda_L=\{0,\ldots,L-1\}^d$ moving in the field $\a$.  
At most one particle is
allowed at each site.  A particle at $x\in\Lambda $ attempts to 
jump to nearest neighbor site $y\in\Lambda$ at rate 
$$1+e^{\alpha_y-\alpha_x}.$$  
If there is no particle in the way the particle is
allowed to jump.  
However if there is a particle in the way, the jump is suppressed,
and everything starts again.  
All the particles are doing this independently of each
other, and since time is continuous one can ignore 
the occasion of two particles
trying to jump onto each other simultaneously.

{\bf Bernoulli-Laplace.}  This is introduced as a tool for proving
results about the Kawasaki dynamics.  For any $\Lambda\subset\ZZ^d$,
the Dirichlet form is$$
D_{BL}(\Lambda, N; f)=E_{\mu_{\Lambda,N}}
[ \sum_{
x,y\in\Lambda} (f(T_{x,y} \eta) - f(\eta))^2].
$$
The difference between Bernoulli-Laplace and Kawasaki is that in Kawasaki
only nearest neighbour jumps are allowed, but in Bernoulli-Laplace we allow
jumps to any site.

The following result has been recently been obtained by Caputo [C].

{\bf Theorem 1.} {\sl Let $0\le K<\infty$.  
 There exists a constant $C=C(K)<\infty$ such that
for any field $\a$ with $-K\le \a_x\le K$, 
 any $\Lambda$, and any $0\le N\le |\Lambda|$,
for any $f:\{0,1\}^\Lambda \to \RR$,
$$
{\rm Var}_{\mu_{\Lambda,N}}(f) \le
{C\over |\Lambda|} D_{BL}(\Lambda,N; f).
$$
}

The main result of this article is

{\bf Lemma 1. (Moving Particles Lemma) [QY]}  {\sl
Suppose $\mu$ is an inhomogenous  Bernoulli measure on $\ZZ^1$
 with external field $\alpha$
taking values in $[-K,+K]$, conditioned to have $N$ particles. 
Then 
$$
E_\mu[ (f(T_{1L}\eta) - f(\eta) )^2 ]
\leq e^{13K} L
\sum_{1 \leq x \leq L-1 } E_\mu[ (f(T_{x,x+1} \eta) -f( \eta))^2 ] 
$$
}

The proof is given in section 3.

Once we have the Moving Particles Lemma the spectral gap and two block estimate follow using standard arguments [GPV], [KL], [Q1]. 
  Let $m^K_{x}= Av_{|y-x|\le K} \eta_y$ be the 
empirical particle density in a box of radius $K$ around $x$.  The two block estimate
says that, suitably averaged, such a quantity is not substantially different if measured
on a large microscale, or a small macroscale.  Let $\Lambda_L$ be a cube of side length $L$
with periodic boundary conditions and let ${\cal P}_{\Lambda_L, N}$ be the set of probability
densities with respect to $\mu_{\Lambda_L, N}$.

{\bf Theorem 2. (Two block estimate) }   {\sl 
Let $F$ be a continuous function on $[0,1]$ and $\varphi$ a smooth function on the $d-$dimensional
unit torus.  If $K\to\infty$ as $L\to \infty$ with $K\le \delta L$,
$$
\limsup_{\delta\to 0}\limsup_{L\to \infty} \sup_{0\le N\le L^d\atop
f\in{\cal P}_{\Lambda_L, N} } \Big\{
L^{-d} E_{\mu_{\Lambda_L, n}}[ |\sum_x \varphi(x/N) ( F( m_x^K) - F(m_x^{\delta N})) | f]
$$$$
\qquad\qquad\qquad\qquad\qquad- L^{2-d} D_{Kaw}(\Lambda_L, N; \sqrt{f}) \Big\} \le 0.
$$
}

By a box $\Lambda_L$ of side length $L$ we mean a set
of the form $\{ x\in \ZZ^d : x_i-y_i\in \{0,1,\ldots,L-1\}, ~i=1,\ldots,d\}$
for some $y=(y_1,\ldots, y_d)\in \ZZ^d$.  The following theorem gives the 
spectral gap of the Kawasaki dynamics to correct order.  It is stated in [C] based on the Moving Particles
Lemma above from [QY].

{\bf Theorem 3.}   {\sl  For each $K>0$ there exists a $C<\infty$ such that
for all $\a$ with $-K\le \a_x\le K$, all boxes $\Lambda_L$ of side
length $L$, all $0\le N\le L^d$, and all $f:\{0,1\}^{\Lambda_L}\to \RR$,
$$
E_{\mu_{\Lambda_L,N}} [ ( f-E_{\mu_{\Lambda,N}} [  f])^2]
\le C L^2 D_{Kaw}(\Lambda_L,N; f).
$$
}

{\it Proof.} By Theorem 1 we have$$
E_{\mu_{\Lambda_L,N}}[(f -E_{\mu_{\Lambda_L,N}}[f])^2]  \le
{C\over |\Lambda_L|} E_{\mu_{\Lambda_L,N}}
[ \sum_{x,y\in\Lambda_L} (f(T_{x,y} \eta) - f(\eta))^2]
$$
 For each $x,y\in \Lambda_L$  choose a canonical path
$x=x_1,x_2,\ldots, x_n=y$ with $x_i$ and  $x_{i+1}$ by moving first
in the first coordinate direction, then in the second coordinate 
direction, etc. 
By Lemma 1, we have
$$
E_{\mu_{\Lambda_L, N}}[ 
 (f(T_{xy}\eta) - f(\eta) )^2 ]
\leq e^{13K} n
\sum_{1 \leq i \leq n-1 } E_\mu[ (f(T_{x_i,x_{i+1}} \eta) -f( \eta))^2 ]. 
$$
Summing over $x$ and $y$, noting that $n\le dL$ and that each nearest 
neighbour pair is used for the path between $d(L/2)^{d+1}$ pairs $x$ and $y$
we obtain the result.\hfill$\qed$

\bigskip
{\bf 3.  Proof of the main result
}
\medskip

{\bf  Lemma 2.}   {\sl
 Suppose $\mu$ is a homogeneous
 Bernoulli measure ($\a_x\equiv 0$) on $\ZZ^1$ conditioned to have
$N$ particles.  Let
$k$ be a positive integer and $\rho_x$
a  sequence of positive numbers with $\sum_{x=}^{k-1} \rho_x =1$.
Then 
$$
E_{\mu} [f(T_{1,k}\eta) - f(\eta) ]^2 
\leq  
 \sum_{x =1}^{k-1 } \rho_x^{-1} E_\mu[(f(T_{x,x+1} \eta) -f(\eta))^2].
$$
}

{\bf Proof.} 
By definition,
$
T_{1,k}\eta = T_{1,2}...T_{k-2, k-1} T_{k, k-1}...T_{3,2}T_{2,1} \eta.
$
For notational convenience, let 
$
T_{k+s-1,k+s} = T_{k-s, k-s-1}$ and $\eta_{k+s}=\eta_{k-s}$, $1\leq s \leq k-1$.
We have 
$$
\eta_1(1-\eta_k) \bigr[f(T_{1,k} \eta) -f(\eta) \bigr]
=\sum _{s=1}^{2k-1} \eta_1(1-\eta_k) q_s(\eta)  \bigr [ 
f(T_{s+1, s}...T_{3,2}T_{2,1} \eta)
-f(T_{s, s-1}...T_{2,1} \eta) \bigr ]
$$
where
$$
q_s(\eta) = \cases{1-\eta_{s+1},  & if $s\leq k-1$; \cr 
                   \eta_{2k-s-2},&  if $s \geq k$. }
$$
Note that the factor $q_s(\eta) $ can be added free of charge because
 the summand vanishes exactly when it does.  Also, for any $ 1\leq x\leq k-1$ fixed,
$ q_x + q_{2k-x-2} = 1$. Hence if one lets $\rho_{k+s} = \rho_{k-s-2}$ then,
by Schwarz's ineqality,
$$
\eqalign{
&   \eta_1(1-\eta_k) \bigr [ f( \eta^{1k}) -f(\eta) 
\bigr ]^2  \cr
& \leq
  \eta_1(1-\eta_k)  \sum_{s=1}^{2k-1} \rho_{s-1}^{-1}   
q_s(\eta) \bigr [ f(T_{s+1 s}\cdots T_{21} \eta)
-f(T_{s s-1}\cdots T_{21} \eta) \bigr ]^2  
  \sum_{s=1}^{2k-2} \rho_{s-1} q_s(\eta)   \cr
& =
 \sum_{s=1}^{2k-2} \rho_{s-1}^{-1}  \eta_1(1-\eta_k) q_s(\eta)
\left[ f(T_{s+1 s}\cdots T_{21} \eta)
-f(T_{s s-1}\cdots T_{21} \eta) \right ]^2     
}
$$
Taking expectation and  change variables $ T_{s s-1}\cdots T_{21} \eta \to \eta$ and also
change the index $ s \to 2k-s-2 $ for $s \geq k$ we find that $E^\mu[ \eta_1(1-\eta_k) \bigr [ f( \eta^{1k}) -f(\eta) 
\bigr ]^2]$
 is bounded by the expectation of
$$ \eqalign{ &\sum_{s=1}^{k-1} \rho_{s-1}^{-1}  (1-\eta_k) (1-\eta_{s+1})
\left[ f(T_{s+1, s} \eta)
-f( \eta) \right]^2  +
\sum_{s=1}^{k-2} \rho_{s-1}^{-1}  \eta_k (1- \eta_{s+1})
\left[ f(T_{s+1, s} \eta)
-f( \eta) \right]^2\cr  
&=
 \sum_{x=1}^{k-1} \rho_{x-1}^{-1}   (1-\eta_{x+1})
\left[ f(T_{x+1, x} \eta)
-f( \eta) \right]^2    
}
$$
We have thus proved that 
$$
E_\mu[  \eta_1(1-\eta_k) ( f( T_{1,k}\eta) -f(\eta) 
)^2]  
 \leq
\sum_{x=1}^{k-1} \rho_{x-1}^{-1} E_\mu[  (1-\eta_{x+1})
( f(T_{x+1, x} \eta)
-f( \eta) )^2 ]  
$$
Similiarly, if one uses the particle-hole duality,
$$E_\mu[  \eta_k(1-\eta_1) ( f( T_{1,k}\eta) -f(\eta) 
)^2]  
 \leq
\sum_{x=1}^{k-1} \rho_{x-1}^{-1} E_\mu[  \eta_{x+1}
( f(T_{x+1, x} \eta)
-f( \eta) )^2 ]
$$
The lemma is obtained by adding these two bounds. 

{\bf Proof of lemma 1.} 
Suppose we change the measure $\mu$ to  a new measure 
$\tilde \mu$ by changing  each $\alpha_i$ to the nearest
value of the form $Kj/L$, $j$ an integer.  Since there is always such a point
with $|\alpha_i - Kj/L | \le K(2L)^{-1}$ the Radon-Nikodym derivative $d\mu/d\tilde\mu$
is bounded above and below uniformly by $e^{K/2}$ and $e^{-K/2}$ respectively.  Therefore
at the cost of a factor of $e^K$  we may assume that $\alpha$ takes values
in $\{ Kj/L: j =-L,\ldots,L\}$.
By the same reasoning, at the  price of a factor $e^{4K}$ 
we may assume that $\alpha_0 =\alpha_L~=K$.

Let A be the set 
$$
A= \{ x_i : \alpha_{x_i} = K, i=1,\ldots,k. \}
$$
By definition,
$$
T_{1,L}\eta = 
T_{x_1,x_2}\cdots T_{x_{k-2},x_{ k-1}} T_{x_{k},x_{k-1}}\cdots T_{x_3,x_2}T_{x_2,x_1} \eta.
$$
By Lemma 2 with  $\rho_s^{-1} = { L\over {x_{s+1}-x_{s}}}$,
$$ 
E_\mu[(f(T_{1,L}\eta) - f(\eta))^2]  \leq e^{4K}
\sum_{s=1}^{k-1} { L\over {x_{s+1}-x_{s}}} E_\mu[ 
(f(T_{x_s, x_{s+1}} \eta)
-f( \eta))^2]
$$
We have to bound 
$$
 { L\over {x_{s+1}-x_{s}}} E_\mu[
( f(T_{x_{s} ,x_{s+1}} \eta)
-f( \eta))^2 ]
$$
We are now in the same situation as before except no $\alpha_x$
can take value K when $x_{s} < x < x_{s+1}$. 
Let us change $\alpha_{x_{s}}$
and $\alpha_{x_{s+1}}$ to the value $ K(L-1)/L$. 
The price we pay is a factor
$\exp\{2KL^{-1}\}$. Continuing this procedure we have 
a proof of the lemma.
\hfill$\qed$

\bigskip

\centerline{\bf References}
\smallskip
\frenchspacing

\item{[AHL]} Ambegaokar, V., Halperin, B.I., Langer, J.S.,
 Hopping conductivity in disordered systems, Phys. Rev. B 4 (1971) 2612.

\item{[BE]}
Brak, R., Elliott, R. J., Correlated random walks with random hopping rates, Journal of
Physics - Condensed Matter, 1989 Dec 25, V1 N51: 10299-10319.
, Correlated tracer diffusion in a disordered medium, 
Materials Science and Engineering B - Solid State Materials for Advanced Technology,
1989 Jul, V3 N1-2: 159-162.

\item{[C]} Caputo, P., Spectral gap inequalities in product spaces with conservation laws, preprint.
www.math.tu-berlin.de/stoch/Kolleg/Homepages/caputo/research/workit.html.

\item{[F]} Fritz, J.,
 Hydrodynamics in a symmetric random medium, 
Comm. Math. Phys. {\bf 125}, (1989)
13-25.

\item{[FM]}  Faggionato, A.,  Martinelli, F., Hydrodynamic limit of a disordered lattice gas. ArXiv. math.PR/0302123
\item{[GP]}
 Gartner, P. and Pitis, R.,  Occupancy-correlation corrections in hopping, Phys. Rev. B 45 (1992).
\item{[GPV]} Guo, M. Z., Papanicolaou, G. C., Varadhan, S. R. S., 
 Nonlinear diffusion limit for a system with nearest neighbor interactions, Comm. Math. Phys. 
{\bf 118} (1988) 31-53.

\item{[KPW]}
 Kehr, K. W., Paetzold, O. and Wichmann, T., Collective diffusion
of lattice gases on linear chain with site-energy disorder,  Phys. Lett. A. ??

\item{[K]} Kirkpatrick, S.,
 Classical transport in disordered media:  Scaling and effective-medium theories,
Phys. Rev. Lett. 27 (1971) 1722.

\item{[LY]} Lu, S. L., Yau, H. T., 
 Spectral gap and logarithmic Sobolev inequality for Kawasaki and
Glauber dynamics, 
Comm. Math. Phys. {\bf 156} (1993) 399-433.

\item{[MA]} Miller, A., Abrahams, E.,
 Impurity conduction at low concentrations, Phys. Rev. 120 (1960) 745.

\item{[Q2]}   Quastel, J. Diffusion in disordered media.  Nonlinear stochastic PDEs (Minneapolis, MN, 1994),  65--79, IMA Vol. Math. Appl., 77, Springer, New York, 1996.

\item{[Q1]} Quastel, J.,  Diffusion of color in the simple exclusion process, Comm. Pure Appl. Math.
XLV, 623-679 (1992).

\item{[QY]} Quastel, J., Yau, H.T.  Bulk diffusion in a system with site disorder, manuscript.

\item{[R]} Richards, P. M.,  Theory of one-dimensional hopping conductivity and
diffusion, Phys. Rev. B 16 (1977) 1393-1409.
\item{[Se]} Seppalainen, T.,  Recent results and open problems on the hydrodynamics of disordered asymmetric exclusion and zero-range processes. II Brazilian School of Probability (Portuguese) (Barra de Sah\'y, 1998).  Resenhas  4  (1999),  no. 1, 1--15. 

\item{[Sp]} Spohn, H., Large scale dynamics of interacting particles, Springer-Verlag (1991).

\vfill
\eject

\end